\documentclass[12pt]{article}
\usepackage{amscd}
\begin{document}
\renewcommand{\thesubsection}{\arabic{subsection}}
\newenvironment{proof}{{\bf Proof}:}{\vskip 5mm }
\newenvironment{rem}{{\bf Remark}:}{\vskip 5mm }
\newenvironment{remarks}{{\bf Remarks}:\begin{enumerate}}{\end{enumerate}}
\newenvironment{examples}{{\bf Examples}:\begin{enumerate}}{\end{enumerate}}  
\newtheorem{proposition}{Proposition}[subsection]
\newtheorem{lemma}[proposition]{Lemma}
\newtheorem{definition}[proposition]{Definition}
\newtheorem{theorem}[proposition]{Theorem}
\newtheorem{cor}[proposition]{Corollary}
\newtheorem{conjecture}{Conjecture}
\newtheorem{pretheorem}[proposition]{Pretheorem}
\newtheorem{hypothesis}[proposition]{Hypothesis}
\newtheorem{example}[proposition]{Example}
\newtheorem{remark}[proposition]{Remark}
\newcommand{\llabel}[1]{\label{#1}}
\newcommand{\comment}[1]{}
\newcommand{\sr}{\rightarrow}
\newcommand{\dw}{\downarrow}
\newcommand{\bdl}{\bar{\Delta}}
\newcommand{\zz}{{\bf Z\rm}}
\newcommand{\zq}{{\bf Z}_{qfh}}
\newcommand{\nn}{{\bf N\rm}}
\newcommand{\qq}{{\bf Q\rm}}
\newcommand{\nq}{{\bf N}_{qfh}}
\newcommand{\oo}{\otimes}
\newcommand{\uu}{\underline}
\newcommand{\ih}{\uu{Hom}}
\newcommand{\af}{{\bf A}^1}
\newcommand{\dsr}{\stackrel{\sr}{\scriptstyle\sr}}

\begin{center}
{\Large\bf Homotopy theory of simplicial sheaves in completely
decomposable topologies}\\ 
\vskip 4mm
\vskip 3mm

 {\bf Vladimir Voevodsky}\footnote{Supported by the NSF
grants DMS-97-29992 and DMS-9901219, Sloan Research Fellowship and Veblen
Fund}$^,$\footnote{School of Mathematics, Institute for Advanced Study,
Princeton NJ, USA. e-mail: vladimir@ias.edu}\\
{August 2000}
\end{center}
\vskip 4mm
\tableofcontents
\subsection{Introduction}
There are two substantially different approaches to the homotopy
theory of simplicial (pre-)sheaves. The first one, introduced by Joyal
and developed in full detail by Jardine in \cite{J2} and \cite{J3},
provides us with a proper simplicial model structure on the category
of simplicial sheaves and presheaves and allows to give definitions of
all the important objects of homotopical algebra in the context of
sheaves on all sites. The drawback of this approach, which is a
necessary consequence of its generality, is that some of the important
classes of objects and morphisms have a very abstract definition. In
particular, the existence of fibrant objects can not be proved without
the use of large cardinals and transfinite arguments even in the case
when the underlying site is reasonably small. Another approach was
introduced in \cite{KSB2} for simplicial sheaves on Noetherian
topological spaces of finite dimension. It is much more explicit and
gives a finitely generated model structure on the category of
simplicial sheaves but does not generalize to arbitrary sites.

The goal of this paper is to define a class of sites for which a
generalized analog of the Brown-Gersten approach works. We prove that
for such sites the class of weak equivalences of simplicial presheaves
can be generated (in the sense of $\Delta$-closed classes considered
in \cite{HH0}) by an explicitly given set of generating weak
equivalences. Using the results of \cite{HH0} we conclude that for
sites with sufficiently good cd-topologies the homotopy category of
simplicial sheaves is equivalent to the localization of
$\Delta^{op}T^{\coprod}$ with respect to the $\bdl$-closure of a very
explicit set of morphisms.

In \cite{HH2} we use this description in the case of the cdh-topology on
schemes over a field to prove a comparison theorem for the motivic
homotopy categories in the Nisnevich and the cdh-topologies.  
\vspace{3mm}

Let $C$ be a small category with an initial object $0$ and $P$ be a
set of commutative squares in $C$. One defines the
cd-topology associated with $P$ as the topology generated by coverings
of two sorts
\begin{enumerate}
\item the empty covering of $0$
\item coverings of the form $\{A\sr X, Y\sr X\}$ where the morphisms
$A\sr X$ and $Y\sr X$ are sides of an element of $P$ of the form
\begin{equation}
\llabel{intrsq}
\begin{CD}
B @>>> Y\\
@VVV @VVV\\
A @>>> X
\end{CD}
\end{equation}
\end{enumerate}
A fundamental example of a
cd-topology is the canonical topology on the category of open subsets
of a Noetherian topological space which is associated with the set of
squares of the form
$$
\begin{CD}
U\cap V@>>>U\\
@VVV @VVV\\
V@>>>U\cup V
\end{CD}
$$
A square (\ref{intrsq}) defines, in the way explained in
\cite{HH0}, an object $K_Q$ of $\Delta^{op}C^{\coprod}$ and a morphism
$K_Q\sr X$. Let $W_P$ be the union of the set of morphisms of this
form with the morphism $0'\sr 0$ where $0'$ is the initial object in
$C^{\coprod}$ and $0$ is the initial object in $C$. Elements of $W_P$
are called generating weak equivalences defined by $P$. Denote by
$W_{proj}$ the class of projective weak equivalences of simplicial
presheaves on $C$ i.e. morphisms $f:X\sr Y$ such that the map of
simplicial sets $X(U)\sr Y(U)$ defined by $f$ is a weak equivalence
for each $U$ in $C$. The main theorem of this paper (Theorem
\ref{tmain}) asserts that, under certain conditions, a morphism of
simplicial presheaves on $C$ is a weak equivalence with respect to the
topology associated with $P$ if and only if it belongs to
$cl_{\bdl}(W_P\cup W_{proj})$. As a corollary, we show that, under the
same conditions, the homotopy category of simplicial presheaves on $C$
with respect to the associated topology is equivalent to the
localization of $\Delta^{op}C^{\coprod}$ with respect to
$cl_{\bdl}(W_P)$.

We distinguish three types of cd-structures - complete, bounded and
regular. One of the properties of complete cd-structures is that any
presheaf which takes distinguished squares to pull-back squares is a
sheaf in the associated topology. One of the properties of regular
cd-structures is that any sheaf in the associated topology takes 
distinguished squares to  pull-back squares. Bounded cd-structures are
intuitively the ones where every object has a finite dimension. In
particular any object has finite cohomological dimension with respect
to the associated topology.

For any complete bounded cd-structure we prove an analog of the
Brown-Gersten theorem \cite[Th. 1']{KSB2}. Using this result we show
that for a complete, regular and bounded cd-structure the class of
weak equivalences with respect to the associated topology is generated
by $W_P$. We also define in this case a closed model structure on the
category of simplicial $t_P$-sheaves which is a general version of the
closed model structure introduced in \cite{KSB2}. This closed model
structure is always finitely generated.

In Section \ref{sec10} we consider cd-structures on
categories with fiber products. For any morphism $f:X\sr Y$
in $C$ denote by $\check{C}(f)$ the simplicial object with terms 
$$\check{C}(f)_{n}=X^{n+1}_Y$$ and faces and degeneracy morphisms
given by partial projections and diagonals respectively. The
projections $X^{i}_Y\sr Y$ define a morphism $\check{C}(f)\sr Y$ which
we denote by $\eta(f)$. We show that for a regular cd-structure such
that a pull-back of a distinguished square is a distinguished square
the $\bdl$-closure of the class $W_P$ coincides with the
$\bdl$-closure of the class of morphisms of the form $\eta(\pi)$ for
$t_P$-coverings $\pi$. This result is important for applications where
it is easier to compute a functor on morphisms of the form $\eta(f)$
than on morphisms of the form $p_Q$ (for example if the functor
commutes with fiber products but not with coproducts).

This paper was written while I was a member of the Institute for
Advanced Study in Princeton and, part of the time, an employee of the
Clay Mathematics Institute. I am very grateful to both institutions
for their support. I would also like to thank Charles Weibel who
pointed out a number of places in the previous version of the paper
which required corrections.

\subsection{Topologies defined by cd-structures}
\begin{definition}
\llabel{d1} Let $C$ be a category with an initial object. A 
cd-structure on $C$ is a collection  $P$ of commutative
squares of the form
\begin{equation}
\llabel{eq1}
\begin{CD}
B@>>>Y\\
@VVV @VVpV\\
A @>e>> X
\end{CD}
\end{equation}
such that if $Q\in P$ and $Q'$ is isomorphic to $Q$ then $Q'$ is in $P$.
\end{definition}
The squares of the collection $P$ are called distinguished squares of
${ P}$. One can combine different cd-structures and/or restrict them to
subcategories considering only the
squares which lie in the corresponding subcategory. Note also that if
$C$ is a category with a  cd-structure $P$ and $X$ is an object of
$C$ then $P$ defines a  cd-structure on $C/X$.

We define the topology $t_P$ associated to a cd-structure $P$ as the
smallest Grothendieck topology (see \cite[III.2]{MM}) such that for a
distinguished square of the form (\ref{eq1}) the sieve $(p,e)$
generated by the morphisms
\begin{equation}
\llabel{elcov}
\{p:Y\sr X, e:A\sr X\}
\end{equation}
is a covering sieve and such that the empty sieve is a covering sieve
of the initial object $\emptyset$. The last condition implies that for
any $t_{P}$-sheaf $F$ one has $F(\emptyset)=pt$.  
We define simple coverings as the ones which can be obtained by iterating
elementary coverings of the form (\ref{elcov}). More precisely one
has.
\begin{definition}
\llabel{simplcov} The class $S_{P}$ of simple coverings is the
smallest class of families of morphisms of the form $\{U_{i}\sr
X\}_{i\in I}$ satisfying the following two conditions:
\begin{enumerate}
\item any isomorphism $\{f\}$ is in $S_{P}$
\item for a distinguished square $Q$ of the form (\ref{eq1}) and
families $\{p_i:Y_i\sr
Y\}_{i\in I}$, $\{q_j:A_j\sr A\}_{j\in J}$ in $S_P$ the family $\{p\circ
p_i, e\circ q_j\}_{i\in I, j\in J}$ is in $S_P$. 
\end{enumerate}
\end{definition}
\begin{definition}
\llabel{complete}
A  cd-structure is called complete if any covering sieve of an
object $X$ which is not isomorphic to $\emptyset$ contains
a sieve generated by a simple covering.
\end{definition}
\begin{lemma}
\llabel{l8.5.1}
A cd-structure $P$ on a category $C$ is complete if and only if the
following two conditions hold:
\begin{enumerate}
\item any morphism with values in $\emptyset$ is an isomorphism 
\item for any distinguished
square of the form (\ref{eq1}) and any morphism $f:X'\sr X$ the sieve
$f^*(e,p)$ contains the sieve generated by a simple covering 
\end{enumerate}
\end{lemma}
\begin{proof}
The first condition is necessary because if there is a morphism $U\sr
0$ then the empty sieve on $U$ is the pull-back of a covering sieve
and therefore a covering sieve. Since the empty sieve does not contain
the sieve generated by any simple covering this contradicts the
completeness assumption. The second condition is necessary for obvious
reasons. To prove that these two conditions are sufficient we have to
show that they imply that the class of sieves which contain the sieves
generated by simple coverings satisfy the Grothendieck topology axioms
(see e.g. \cite[Def. 1 p.110]{MM}). The first and the third axioms
hold for the class of simple coverings in any cd-structure. An 
inductive argument shows that  that the conditions of the lemma imply
that the second, stability, axiom holds.
\end{proof}
Lemma \ref{l8.5.1} immediately implies the following three lemmas.
\begin{lemma}
\llabel{completepb} Let $P$ be a  cd-structure such that 
\begin{enumerate}
\item any morphism with values in $\emptyset$ is an isomorphism 
\item for any
distinguished square $Q$ of the form (\ref{eq1}) and any
morphism $X'\sr X$ the square $Q'=Q\times_X X'$ is defined and belongs
to $P$.
\end{enumerate}
Then $P$ is complete.
\end{lemma}
\begin{lemma}
\llabel{l8.5.2} Let $C$ be a category and $P_1$, $P_2$ two complete
cd-structures on $C$. Then $P_1\cup P_2$ is complete.
\end{lemma}
\begin{lemma}
\llabel{l8.5.3}
Let $C$ be a category and $P$ a complete cd-structure on $C$. Then for
any object $U$ of $C$ the cd-structure $P/U$ on $C/U$ defined by $P$
is complete.
\end{lemma}
All simple coverings are finite and therefore the topology associated
to any complete cd-structure is necessarily Noetherian. In
particular, we have the following compactness result.  For an object
$X$ of $C$ denote by $\rho(X)$ the $t_P$-sheaf associated to the
presheaf represented by $X$.
\begin{lemma}
\llabel{compactness} Let $P$ be a complete cd-structure. Then for any
$U$ in $C$ the sheaf $\rho(U)$ is a compact object of $Shv(C,t_P)$
i.e. for any filtered system of sheaves
$F_{\alpha}$ one has
\begin{equation}
\llabel{eq8.2.0}
Hom(\rho(U), colim(F_{\alpha}))=colim Hom (\rho(U), F_{\alpha})
\end{equation}
\end{lemma}
\begin{proof}
By definition of associated sheaf and the Yoneda Lemma the right hand
side of (\ref{eq8.2.0}) is isomorphic to $colim(F_{\alpha}(U))$ and
the left hand side to $(colim F_{\alpha})(U)$. Therefore, to prove the
lemma it is enough to show that the colimit of the family $F_{\alpha}$
in the category of presheaves coincides with its colimit in the
category of sheaves i.e. that the colimit in the category of
presheaves is a sheaf. Since every covering
has a finite refinement this follows from the fact that finite limits
in the category of sets commute with filtered colimits.
\end{proof}
\begin{lemma}
\llabel{isasheaf} Let $P$ be a complete  cd-structure and $F$ a
presheaf on $C$ such that $F(\emptyset)=*$ and for any distinguished
square $Q$ of the form (\ref{eq1}) the square 
\begin{equation}
\llabel{FQ}
F(Q)=\left(
\begin{CD}
F(X)@>>>F(Y)\\
@VVV @VVV\\
F(A)@>>>F(B)
\end{CD}
\right)
\end{equation}
is pull-back. Then $F$ is a sheaf in the associated topology.
\end{lemma}
\begin{proof}
Let us say that a section of $F$ on a sieve $J=\{p_U:U\sr X\}$ is a
collection of sections $a_U\in F(U)$ satisfying the condition
$F(g)(a_U)=a_V$ for any morphism $g:V\sr U$ over $X$.  A presheaf $F$
is a sheaf if for any object $X$ of $C$, any covering sieve
$J=\{p_U:U\sr X\}$ of $X$ and any section $(a_U)$ of $F$ on $J$ there
exists a unique section $a_X\in F(X)$ such that $a_U=F(p_U)(a_X)$ for
all $p_U\in J$. Observe first that our condition on $F$ implies that
for any simple covering $\{U_i\sr X\}$ the map $F(X)\sr \prod F(U_i)$
is injective. Together with the completeness condition it implies
that if $a_X$ exists it is unique i.e. $F$ is a separated
presheaf. The same condition implies that it is sufficient to prove
existence for sieves of the form $(p_i)$ where $\{p_i:U_i\sr X\}$ is a
simple covering. Consider the class $S_F$ of simple coverings such
that for $\{p_i\}\in S_F$ the existence condition holds for any
section $(a_{U})$. Let us show that it coincides with
$S_P$. It clearly contains isomorphisms.  Let $Q$, $p_i$ and $q_j$ be as in
Definition \ref{simplcov}(2) and assume that $(p_i), (q_j)\in
S_F$. Our section restricted to $(p_i)$ and $(q_j)$ defines two section
$a_Y\in F(Y)$ and $a_A\in F(A)$. If $u:W\sr Y$, $v:W\sr A$ are
morphisms such that $pu=ev$ then there exists a covering sieve $J$ on
$W$ contained in both $u^*((p_i))$ and $v^*((q_j))$. The sections
$F(u)(a_Y)$ and $F(v)(a_A)$ coincide on elements from $J$ and since
$F$ is separated they coincide in $F(W)$. We conclude that $a_Y$ and
$a_A$ coincide on $B$ and thus our condition on $F$ implies that there
exists $a_X$ such that $a_A=F(e)(a_X)$ and $a_Y=F(p)(a_X)$.
\end{proof}
\begin{definition}
\llabel{cdstr}
A cd-structure $P$ is called regular if for any $Q\in P$ of the form
(\ref{eq1}) one has
\begin{enumerate}
\item $Q$ is a pull-back square
\item $e$ is a monomorphism
\item the morphism of
sheaves 
\begin{equation}
\llabel{eq2new}
\Delta\coprod \rho(e_B)^2:\rho(Y)\coprod
\rho(B)\times_{\rho(A)}\rho(B)\sr \rho(Y)\times_{\rho(X)}\rho(Y)
\end{equation}
is surjective.
\end{enumerate}
\end{definition}
The following three lemmas are straightforward.
\begin{lemma}
\llabel{cdstrpb}
Let $P$ be a cd-structure such that for any distinguished square
$Q$ of the form (\ref{eq1}) in $P$ one has
\begin{enumerate}
\item $Q$ is a pull-back square
\item $e$ is a monomorphism
\item the objects $Y\times_X Y$ and
$B\times_A B$ exist and the derived square
\begin{equation}
\llabel{eq2}
d(Q)=\left(
\begin{CD}
B@>>>Y\\
@VVV @VVV\\
B\times_A B@>>>Y\times_X Y
\end{CD}
\right)
\end{equation}
where the vertical arrows are the diagonals, is distinguished.
\end{enumerate}
Then $P$ is a regular cd-structure.
\end{lemma}
\begin{lemma}
\llabel{l8.5.4}
Let $C$ be a category and $P_1$, $P_2$ two regular
cd-structures on $C$. Then $P_1\cup P_2$ is regular.
\end{lemma}
\begin{lemma}
\llabel{l8.5.5}
Let $C$ be a category and $P$ a regular cd-structure on $C$. Then for
any object $U$ of $C$ the cd-structure $P/U$ on $C/U$ defined by $P$
is regular.
\end{lemma}
\begin{example}\rm
\llabel{toy0}
Define the toy cd-structure as follows. The category $C$ is the
category corresponding to the diagram
\begin{equation}
\llabel{toysq}
\begin{CD}
B@>>>Y\\
@VVV @VVV\\
A @>e>> X
\end{CD}
\end{equation}
together with an additional object $\emptyset$ which is an initial
object and such that any morphism with values in  $\emptyset$ is
identity. The distinguished squares are the square (\ref{toysq}) and
the squares 
\begin{equation}
\llabel{toysq2}
\begin{CD}
\emptyset @>>> \emptyset\\
@VVV @VVV\\
A@>>>A
\end{CD}
\,\,\,\,\,\,
\begin{CD}
\emptyset@>>>\emptyset\\
@VVV @VVV\\
\emptyset@>>>\emptyset
\end{CD}
\end{equation}
One verifies easily that the toy cd-structure is complete 
and regular.
\end{example}
\begin{proposition}
\llabel{p8.2.2}
For any regular cd-structure $P$, distinguished square $Q$ of the form
(\ref{eq1}) and a sheaf $F$ in the associated topology the square
of sets $F(Q)$ is a pull-back square.
\end{proposition}
\begin{proof}
The map $F(X)\sr F(A)\times F(Y)$ is
a monomorphism since $\{Y\sr X, A\sr X\}$ is a covering and its image
coincides with the equalizer of the maps 
$$Hom(\rho(Y)\coprod \rho(A),F)\sr Hom((\rho(Y)\coprod
\rho(A))\times_{\rho(X)}(\rho(Y)\coprod \rho(A)), F)$$
defined by the projections. Since (\ref{eq2new}) is an epimorphism and
$e$ is a monomorphism this equalizer coincides with the equalizer of
the maps $F(Y)\times F(A)\sr F(B)$ defined by the morphisms from $B$
to $A$ and $Y$.
\end{proof}
Proposition \ref{p8.2.2} has the following immediate corollaries.
\begin{cor}
\llabel{l1}
For any regular cd-structure and any distinguished square $Q$ of the form
(\ref{eq1}) the square
\begin{equation}
\llabel{eq3}
\rho(Q)=\left(
\begin{CD}
\rho(B)@>>>\rho(Y)\\
@VVV @VVV\\
\rho(A)@>>>\rho(X)
\end{CD}
\right)
\end{equation}
is a push-out square.
\end{cor}
\begin{cor}
\llabel{regcompl} Let $P$ be a complete regular cd-structure. Then a
presheaf $F$ is a sheaf in the associated topology if and only if
$F(\emptyset)=pt$ and for any distinguished square $Q$ of the form
(\ref{eq1}) the square (\ref{FQ}) is pull-back.
\end{cor}
For a sheaf $F$ let $\zz(F)$ be the sheaf of abelian groups freely
generated by $F$.
\begin{lemma}
\llabel{MV1}
Let $P$ be a regular cd-structure and $Q$ be a distinguished square of
the form (\ref{eq1}). Then sequence of sheaves of abelian groups
\begin{equation}
\llabel{seq}
0\sr \zz(\rho(B))\sr \zz(\rho(A))\oplus \zz(\rho(Y))\sr
\zz(\rho(X))\sr 0
\end{equation}
is exact.
\end{lemma}
\begin{proof}
For any site the functor $F\mapsto \zz(F)$ takes colimits to colimits
and monomorphisms to monomorphisms. Since $P$ is regular, the morphism
$B\sr Y$ is a monomorphism and therefore the sequence (\ref{seq}) is
exact in the first term. The quotient $\zz(\rho(A))\oplus
\zz(\rho(Y))/\zz(\rho(B))$ is the colimit of the diagram
$$
\begin{CD}
\zz(\rho(B))@>>>\zz(\rho(Y))\\
@VVV\\
\zz(\rho(A))
\end{CD}
$$
and since $\zz(-)$ is right exact, Corollary \ref{l1} implies that it is
$\zz(\rho(X))$.
\end{proof}
\begin{lemma}
\llabel{MV2}
Let $P$ be a regular cd-structure and $F$ a $t_P$-sheaf of abelian
groups on $C$. Then there is a function which takes any distinguished
square of the form 
(\ref{eq1}) to a family of homomorphisms 
$$\partial_Q:H^i(B,F)\sr H^{i+1}(X,F)$$
and which satisfies the following conditions
\begin{enumerate}
\item for a map of squares $Q'\sr Q$ the square 
\begin{equation}
\begin{CD}
H^i(B,F) @>\partial_Q>> H^{i+1}(X,F)\\
@VVV @VVV\\
H^i(B',F) @>\partial_{Q'}>> H^{i+1}(X',F)
\end{CD}
\end{equation}
commutes
\item the sequence of abelian groups
\begin{equation}  
H^i(X,F)\sr H^i(A,F)\oplus H^i(Y,F)\sr H^i(B,F)\sr
H^{i+1}(X,F)
\end{equation}
is exact.
\end{enumerate}
\end{lemma}
\begin{proof}
Our result follows from Lemma \ref{MV1} and the fact that for any
site $T$, sheaf $F$ on $T$ and an object $X$ of $T$ one has
$$H^n(X,F)=Ext^n(\zz(\rho(X)),F).$$ 
\end{proof}
Our next goal is to define dimension for objects of a category with a
cd-structure and in particular introduce a class of cd-structures of
(locally) finite dimension. We start with the following auxiliary
definition.
\begin{definition}
\llabel{density} A density structure on a category $C$ with an initial
object is a function which assigns to any object $X$ a sequence
$D_0(X),D_1(X),\dots$ of families of morphisms 
which satisfies the following conditions:
\begin{enumerate}
\item $X$ is the codomain of elements of $D_i(X)$ for all $i$
\item $(\emptyset\sr X)\in D_0(X)$ for all $X$
\item isomorphisms belong to $D_i$ for all $i$
\item $D_{i+1}\subset D_i$
\item if $j:U\sr V$ is in $D_i(V)$ and $j':V\sr X$ is in $D_i(X)$ then
$j'\circ j:U\sr X$ is in $D_i(X)$
\end{enumerate}
\end{definition}
By abuse of notation we will write $U\in D_n(X)$ instead of $(U\sr
X)\in D_n(X)$. A density structure is said to be locally of finite
dimension if for any $X$ there exists $n$ such that any element of
$D_{n+1}(X)$ is an isomorphism. The smallest such $n$ is called the
dimension of $X$ with respect to $D_*(-)$.
\begin{definition}
\llabel{redsq} Let $C$ be a category with a cd-structure $P$ and a
density structure $D_*(-)$. A distinguished square $Q$ of the form
(\ref{eq1}) is called reducing (with respect to $D_*$) if for any
$i\ge 0$, and any $B_0\in D_i(B)$, $A_0\in D_{i+1}(A)$, $Y_0\in
D_{i+1}(Y)$ there exist $X'\in D_{i+1}(X)$, a distinguished square
\begin{equation}
\llabel{eq1prime}
Q'=\left(
\begin{CD}
B'@>>>Y'\\
@VVV @VVpV\\
A' @>e>> X'
\end{CD}
\right)
\end{equation}
and a morphism $Q'\sr Q$ which coincides with the morphism $X'\sr X$
on the lower right corner and whose other respective components factor
through $B_0, Y_0, 
A_0$.
\end{definition}
We say that a distinguished square $Q'$ is a refinement of a
distinguished square $Q$ if there is a morphism $Q'\sr Q$ which is the
identity on the lower right corner.
\begin{definition}
\llabel{reddn} A density structure $D_*(-)$ is said to be a reducing
density structure for a  cd-structure $P$ if any
distinguished square in $P$ has a refinement which is reducing with
respect to $D_*(-)$.  
A cd-structure is called bounded if there exists a reducing density
structure of locally finite dimension for it.
\end{definition}
The class of reducing squares for a given  cd-structure and a
density structure is again a  cd-structure which we call the
associated reducing cd-structure. If the density structure is reducing
the reducing cd-structure generates the same topology as the original
one. In this case one  cd-structure is complete if and only if the
other is. The following two lemmas are straightforward.
\begin{lemma}
\llabel{l8.5.6} Let $C$ be a category and $P_1$, $P_2$ two
cd-structures on $C$ bounded by the same density structure $D$. Then
$P_1\cup P_2$ is bounded by $D$.
\end{lemma}
\begin{lemma}
\llabel{l8.5.7} Let $C$ be a category and $P$ a cd-structure on $C$
bounded by a density structure $D$. Then for any object $U$ of $C$ the
cd-structure $P/U$ on $C/U$ is bounded by the density structure $D/U$.
\end{lemma}
\begin{example}\rm
Let $C$ be the category of Example \ref{toy0}. Define the density
structure on $C$ setting
$$D_n(\emptyset)=Id \mbox{ for $n\ge 0$}$$
$$D_0(A)=\{Id,\emptyset\}\,\,\,\, D_n(A)=\{Id\} \mbox{ for $n>0$}$$
$$D_0(B)=\{Id,\emptyset\}\,\,\,\, D_n(B)=\{Id\} \mbox{ for $n>0$}$$
$$D_0(Y)=\{Id,\emptyset\}\,\,\,\, D_n(Y)=\{Id\} \mbox{ for $n>0$}$$
$$D_0(X)=\{Id, e, \emptyset\}\,\,\,\, D_1(X)=\{Id, e\}\,\,\,\,
D_n(X)=\{Id\} \mbox{ for $n>1$}$$
Every distinguished square of the toy cd-structure is reducing with
respect to this density structure and in particular the density
structure is reducing. Since the dimension of all objects with respect
to $ D$ is $\le 1$ the toy cd-structure is bounded.
\end{example}
\begin{theorem}
\llabel{cohdim}
Let $P$ be a complete regular cd-structure
bounded by a density structure $D$ and $X$ an object of $C$.
Then for any $t_P$-sheaf of abelian groups $F$ on $C/X$ one has
$$H^n_{t_P}(X,F)=0$$
for $n>dim_{D} X$.
\end{theorem}
\begin{proof}
Replacing $C$ by $C/X$ we may assume that $F$ is defined on $C$. We
will show that for any $X$, any $n$ and any class $a\in H^n(X,F)$
there exists an element $j:U\sr X$ of $D_n(X)$ such that
$j^*(a)=0$. We do it by induction on $n$. For $n=0$ the statement
follows from the fact that $\emptyset\sr X$ is in $D_0(X)$ and for any
$t_P$-sheaf $F$ of abelian groups one has $F(\emptyset)=0$.

Replacing $P$ by the class of the reducing squares with respect to
$D$ we may assume that any square in $P$ is reducing.  Let $a$ be
an element in $H^n$ and $n>0$. Since the sheaves associated to the
cohomology presheaves are zero there exists a $t_P$-covering
$\{p_i:U_i\sr X\}$ such that $p_i^*(a)=0$ for all $i$. Since $P$ is
complete this covering has a simple refinement. Let $S$ be the class
of simple coverings such that if $a$ is a class vanishing on an
element of $S$ then it vanishes on an element of $D_n$. It clearly
contains isomorphisms. Thus to show that it coincides with the whole
$S_P$ it is enough to check that it satisfies the condition of
Definition \ref{simplcov}(2). Let $Q$ be a square of the form
(\ref{eq1}) and $\{p_i:Y_i\sr Y\}$, $\{q_j:A_j\sr A\}$ be elements of
$S$. Then there exist monomorphisms $Y_0\sr Y$ and $A_0\sr A$ in
$D_n(Y)$ and $D_n(A)$ respectively such that $a$ restricts to zero on
$A_0$ and $Y_0$. Setting $B_0=B$ and using the definition of a
reducing square we see that there is an element $X'\sr X$ of $D_n(X)$,
a distinguished square $Q'$ based on $X'$ and a morphism $Q'\sr Q$,
which coincides with the embedding $X'\sr X$ on the lower right
corner, such that the restriction $a'$ of $a$ to $X'$ vanishes on $A'$
and $Y'$. By Lemma \ref{MV2} this implies that $a'=\partial_{Q'}(b')$
where $b\in H^{n-1}(B')$. By induction there is an element $B_0\sr B$
in $D_{n-1}(B)$ such that $b'$ vanishes on $B_0$. Applying again the
definition of a reducing square to $Q'$ with respect to $B_0$,
$Y_0=Y'$ and $A_0=A'$ and using the naturality of homomorphisms
$\partial_{Q}$ we conclude that there is an element $X''\sr X'$ in
$D_n(X')$ such that the restriction of $a'$ to $X''$ is zero. By
definition of a density structure the composition $X''\sr X$ is in
$D_n(X)$ which finishes the proof.
\end{proof}

\subsection{Flasque simplicial presheaves  and local equivalences}
\begin{definition}
\llabel{d2} Let $C$ be a category with a cd-structure $P$. A
B.G.-functor on $C$ with respect to $P$ is a family of contravariant
functors $T_q$, $q\ge 0$ from $C$ to the category of pointed sets
together with pointed maps $\partial_Q:T_{q+1}(B)\sr T_{q}(X)$ given
for all distinguished squares of the form (\ref{eq1}) such that the
following two conditions hold:
\begin{enumerate}
\item the morphisms $\partial_Q$ are natural with respect to morphisms
of distinguished squares
\item for any $q\ge 0$ the sequence of pointed sets
$$T_{q+1}(B)\sr T_{q}(X)\sr T_q(A)\times T_q(Y)$$
is exact.
\end{enumerate} 
\end{definition}
The following theorem is an analog of \cite[Th 1']{KSB2}.
\begin{theorem}
\llabel{analog} Let $C$ be a category with a bounded complete 
cd-structure $P$. Then for any
B.G.-functor $(T_q,\partial_Q)$ on $C$ such that the $t_P$-sheaves
associated to $T_q$ are trivial (i.e. isomorphic to the point sheaf
$pt$) and $T_q(\emptyset)=pt$ for all $q$ one has $T_q=pt$ for all $q$.
\end{theorem}
\begin{proof}
Replacing $P$ with the corresponding reducing  cd-structure we may
assume that all distinguished squares of $P$ are reducing.  Let $T_q$
be a B.G.-functor such that the sheaves $aT_q$ associated to $T_q$'s
are trivial. Let us show that for any $d\ge 0$, $q\ge 0$, $X$ and
$a\in T_q(X)$ there exists $j:U\sr X$ in $D_d(X)$ such that
$T_q(j)(a)=*$. We prove it by induction on $d$. For $d=0$ the
statement follows from the fact that $(\emptyset\sr X)\in D_0(X)$ and
$T_q(\emptyset)=*$. Assume that the statement is proved for $d$ and
all $q$ and $X$. Let $a\in T_{q}(X)$ be an element. Then by our
assumption there exists a covering sieve $J$ such that for any $p:U\sr
X$ in $J$ one has $T_q(p)(a)=*$. Since $P$ is
complete $J$ contains a sieve of the form $(p_i)$ for a 
simple covering $\{p_i:U_i\sr X\}$.  Therefore it is sufficient to
show that for any simple covering $(p_i)$ and $a\in T_q(X)$
such that $T_q(p_i)(a)=*$ there exists $(j:U\sr X)\in D_{d+1}(X)$ such
that $T_q(j)(a)=*$.  Let $S$ be the class of simple coverings
$(p_i)$ such that for any $a\in T_q(X)$ such that $T_q(p_i)(a)=*$
there exists $(j:U\sr X)\in D_{d+1}(X)$ such that $T_q(j)(a)=*$.  It
contains isomorphisms since any isomorphism is in $D_{d+1}(X)$ by
definition of a density structure. Let $Q$ be a distinguished
square, $(p_i)$ and $(q_j)$ be as in Definition \ref{simplcov}(2) and
$(p_i)$ and $(q_j)$ are in $S$. Let us show that $(p\circ p_i, e\circ
q_j)$ is in $S$. Given an element $a\in T_q(X)$ such that its
restrictions to $Y_i$ and $A_j$ are trivial we can find $Y_0\in D_{d+1}(Y)$
and $A_0\in D_{d+1}(A)$ such that the restriction of $a$ to $A_0$ and
$Y_0$ is trivial. Set $B_0=B$. Since $Q$ is reducing we can find a map
$Q'\sr Q$ such that $X'\in D_{d+1}(X)$ and the restriction of $a$ to
$Y'$ and $A'$ is trivial. Let $a'$ be the restriction of $a$ to
$X'$. Since $T_q$ is a B.G-functor we have $a'=\partial(b')$ for some
$b'\in T_{q+1}(B')$. By induction there exists $B'_0\in D_{d}(B)$ such
that the restriction of $b'$ to $B'_0$ is trivial. Set $Y'_0=Y'$,
$A'_0=A'$. Since $Q'$ is reducing we can find $Q''\sr Q'$ such that
$X''\in D_{d+1}(X')$ and $B''\sr B'$ factors through $B'_0$. Since
$\partial$ commutes with morphisms of distinguished squares we
conclude that the restriction of $a'$ to $X''$ is trivial. Finally
observe that since $X''\in D_{d+1}(X')$ and $X'\in D_{d+1}(X)$ one has
$X''\in D_{d+1}(X)$ by Definition \ref{density}(4). We conclude that
$S=S_{P}$ which finishes the proof.
\end{proof}
\begin{definition}
\llabel{d4} Let $C$ be a category with a  cd-structure
$P$. A simplicial presheaf ${ F}$ on $C$ is called
flasque with respect to $P$ if ${ F}(\emptyset)$
is contractible and for any
distinguished square of the form (\ref{eq1}) the square of simplicial
sets
\begin{equation}
\llabel{pbs}
{ F}(Q)=\left(
\begin{CD}
{ F}(X)@>>>{F}(Y)\\
@VVV @VVV\\
{F}(A)@>>>{F}(B)
\end{CD}
\right)
\end{equation}
is a homotopy pull-back square.
\end{definition}
Recall from \cite{HH0} that for a square $Q$ of the form (\ref{eq1})
in a category with finite coproducts we denote by $K_Q$ the simplicial
object given by the elementary exact square
\begin{equation}
\llabel{kq}
\begin{CD}
B\coprod B@>>>B\times\Delta^1\\
@VVV @VVV\\
A\coprod Y@>>>K_Q
\end{CD}
\end{equation}
and by $p_Q:K_Q\sr X$ the obvious morphism.  
\begin{lemma}
\llabel{l8.7.1} Let $F$ be a simplicial presheaf on $C$ such that for
any $U$ in $C$ the simplicial set $F(U)$ is Kan, $F(0)$ is
contractible and for any $Q\in P$ of the form (\ref{eq1}) the map of
simplicial sets
$$F(X)=S(X, F)\sr S(K_Q, F)$$
defined by $p_Q$ is a weak equivalence. Then $F$ is flasque.
\end{lemma}
\begin{proof}
The simplicial set $S(K_Q,F)$ is given by the pull-back square
$$
\begin{CD}
S(K_Q, F) @>>> F(B)^{\Delta^1}\\
@VVV @VVV\\
F(Y)\times F(A) @>>> F(B)\times F(B)
\end{CD}
$$
Therefore, if $F(B)$ is a Kan simplicial set it is a model for the
homotopy limit of the diagram $(F(A)\sr F(B), F(Y)\sr F(B))$ which
implies the statement of the lemma.
\end{proof}
For a presheaf $F$ denote by $aF$ the associated sheaf in the
$t_P$-topology. Recall that a morphism of simplicial presheaves
$f:F\sr G$  is called a weak equivalence with respect to 
the topology $t_P$ if one has:
\begin{enumerate}
\item the morphism $a\pi_0(F)\sr a\pi_0(F)$ defined by
$f$ is an isomorphism
\item for any object $X$ of $C$, any $x\in F(X)$ and any $n\ge
1$ the morphism of associated sheaves $a\pi_n(F,x)\sr
a\pi_n(F,f(x))$ on $C/X$ defined by $f$ is an isomorphism.
\end{enumerate}
Below we call morphisms which are weak equivalences with respect to a
given topology $t$ ``$t$-local weak equivalences'' or, if no confusion
is possible, simply local weak equivalences. Weak equivalences with
respect to the trivial topology are morphisms $f:X\sr Y$ such that for
any $U$ in $C$ the map of simplicial sets $X(U)\sr Y(U)$ defined by
$f$ is a weak equivalence. We call these morphisms projective weak
equivalences.
%
\begin{lemma}
\llabel{key} Let $C$ be a category with a complete bounded
cd-structure $P$. A morphism $f:F\sr G$ of flasque simplicial
presheaves is a $t_P$-local  weak equivalence if and only if it is a
projective weak equivalence.
\end{lemma}
\begin{proof}
Most of the proof is copied from \cite[Lemma 3.1.18]{MoVo}.  The if
part is obvious. Assume that $f$ is a $t_P$-local weak
equivalence. Using the closed model structure on the category of
simplicial presheaves or an appropriate explicit construction we can
find a commutative diagram of simplicial presheaves
\begin{equation}
\begin{CD}
F@>>>G\\
@VVV @VVV\\
F'@>>>G'
\end{CD}
\end{equation}
such that for any $X$ in $C$, the maps 
$F(X)\sr F'(X)$ and  $G(X)\sr G'(X)$
are weak equivalences  of simplicial sets and the map
$F'(X)\sr G'(X)$ is a Kan fibration of
Kan simplicial sets. Replacing $F$, $G$ by $F'$,
$G'$ we may assume that the maps
$F(X)\sr G(X)$ are Kan fibrations between
Kan simplicial sets.

It is sufficient to prove that for any $X$ in 
$C$ and $y\in G(X)$ 
the fiber $K(X)$ of the map $F(X)\sr G(X)$ over $y$ is
contractible ({\it i.e.} weakly equivalent to point
and in particular non empty).
The simplicial presheaf 
$$(p:U\sr X)\mapsto fiber_{G(p)(y)}(F(U)\sr G(U))$$ on $C/X$, clearly
has the B.G.-property with respect to the induced cd-structure on
$C/X$. It is also obvious that a cd-structure on $C/X$ induced by a
complete (resp. bounded) cd-structure is complete
(resp. bounded). Therefore, it is sufficient to prove the lemma for
$G=pt$ in which case we have to show that for any $X$ the (Kan)
simplicial set $F(X)$ is contractible. In addition we may
assume that $C$ has a final object $pt$.

Assume first that $F(pt)\ne\emptyset$ and
let $a\in F_0(pt)$ be an element. Consider the family of functors
$T_q$ on $C$ of the form 
$$X\mapsto \pi_q(F(X),a_{|X}).$$
It is a B.G.-functor and the associated 
$t_P$-sheaves are trivial
since $F\sr pt$ is a weak equivalence. We also have
$T_q(\emptyset)=*$. Thus $T_q(X)=pt$ for all $X$ by Theorem
\ref{analog}. 

It remains to prove that $F(pt)$ is not empty. We already know that
for any $X$ such that $F(X)$ is not empty it is contractible. Since
$P$ is complete the sheaf associated to $\pi_0(F)$ is $pt$ there
exists a simple covering $\{U_i\sr pt\}$ such that $F(U_i)\ne
\emptyset$. Consider the class $S$ of simple coverings $\{U_i\sr X\}$
such that if $F(U_i)\ne \emptyset$ for all $i$ then $F(X)\ne
\emptyset$. Let us show that it coincides with the class of all simple
coverings. It clearly contains isomorphisms. Let $Q$, $p_i$, $q_j$ be
as in Definition \ref{simplcov}(2) and $(p_i),(q_j)$ be in $S$. Then
if $F(U_i)$, $F(V_j)$ are non empty for all $i,j$ then by assumption
$F(Y)$ and $F(A)$ are nonempty. Thus $F(B)$ is nonempty and
therefore all these simplicial sets are contractible. It remains to
note that if in a homotopy pull-back square of the form (\ref{pbs})
$F(A)$, $F(Y)$ and $F(B)$ are contractible then $F(X)$ is non empty.
\end{proof}
Let $W_P^{0}$ be the collection of morphisms of the form $K_Q\sr X$
where $Q\in P$. Denote by $W_P$ the union of $W_P^0$ with the morphism
$\emptyset'\sr \emptyset$ where $\emptyset'$ is the initial object of
$PreShv(C)$ and $\emptyset$ is the presheaf represented by the initial
object of $C$. The elements of $W_P$ are called the generating weak
equivalences defined by $P$.
\begin{lemma}
\llabel{whenis} Consider a commutative square $Q$ of the form
(\ref{eq1}) in the category of presheaves on a site and assume that
$e$ is a monomorphism. Then the square of sheaves associated to $Q$
is a push-out square if and only if the morphism of simplicial
presheaves $K_Q\sr X$ is a local weak equivalence.
\end{lemma}
\begin{proof}
For each $U$ in $C$ the simplicial set $K_Q(U)=K_{Q(U)}$ has the
property that
$$\pi_0(K_{Q}(U))=Y(U)\coprod_{B(U)}A(U)$$
This implies that $\pi_0(K_Q)=Y\coprod_B A$. Therefore, we have
$$a\pi_0(K_Q)=a(Y\coprod_B A)$$
and we conclude that if $K_Q\sr X$ is a local equivalence then $Q$ is
a push-out square. Assume now that $e$ is a monomorphism. Then,
for each $U$ in $C$ the simplicial set $K_{Q(U)}$ is weakly equivalent
to the the (simplicial) set $Y(U)\coprod_{B(U)}A(U)$ i.e. the morphism
$K_{Q}\sr Y\coprod_B A$ is a projective weak equivalence. Since the
associated sheaf functor takes projective weak equivalences to local
weak equivalences we conclude that for a push-out square $Q$ the
morphism $K_Q\sr X$ is a local weak equivalence.
\end{proof}
Combining Lemma \ref{whenis} with Corollary \ref{l1} we get the following
result.
\begin{lemma}
\llabel{es} Let $P$ be a regular cd-structure. Then all elements of
$W_P$ are local weak equivalences for the associated topology.
\end{lemma}
We are going to prove now that for a complete regular and bounded
cd-structure the set $W_P$ generates the class of $t_P$-local weak
equivalences. Let us recall first that in \cite{HH0} we studied two
ways to saturate a given class of morphisms $E$ in a category of the
form $\Delta^{op}R(C)$ where $R(C)$ is the category of radditive
functors on a category $C$ with finite coproducts. The first one is to
consider the class of E-local equivalences $cl_l(E)$. The second is to
consider the $\bdl$-closure $cl_{\bdl}(E\coprod W_{proj})$ of the
union of $E$ with the class of projective weak equivalences.  Since
presheaves on $C$ are radditive functors on the category
$C^{\coprod_{<\infty}}$ obtained from $C$ by adding free finite
coproducts the results of that paper can be applied in our setting.
\begin{proposition}
\llabel{p7.31.5}
For any cd-structure we have
\begin{equation}
\llabel{eq7.31.2}
cl_{l}(W_P)=cl_{\bdl}(W_P\cup W_{proj})
\end{equation}
\end{proposition}
\begin{proof}
By \cite[Remark 4.3.9]{HH0} the category $C^{\coprod_{<\infty}}$
satisfies the condition of \cite[Theorem 4.3.7]{HH0}. Therefore it is
sufficient to show that the domains and codomains of elements of $W_P$
satisfy the conditions of this theorem i.e. that domains and codomains
of elements of $W_P$ are cofibrant in the projective model structure
and reliably compact.  For the morphism $0'\sr 0$ this is
obvious. Consider a morphism of the form $p_Q:K_Q\sr X$ for a
distinguished square $Q$. Since $X$ is a representable presheaf it is
cofibrant and reliably compact. The object $K_Q$ is given by the
push-out square
$$
\begin{CD}
B\coprod B @>>> A\coprod Y\\
@VVV @VVV\\
B\times\Delta^1 @>>> K_Q
\end{CD}
$$
Since $B$ is a representable presheaf it is cofibrant and therefore
the left vertical arrow in this square is a cofibration (by
\cite[Proposition 3.2.11]{HH0}). Since $A$ and $Y$ are cofibrant we
conclude that $K_Q$ is cofibrant.  The same square implies that $K_Q$
is reliably compact.
\end{proof}
\begin{lemma}
\llabel{l7.31.3}
For any site $T$ the class of local weak equivalences in the category
$\Delta^{op}PreShv(T)$ is $\bdl$-closed.
\end{lemma}
\begin{proof}
Observe first of all that if $T$ has sufficiently many points then the
statement of the lemma follows immediately from the corresponding
result for simplicial sets (see e.g. \cite[Lemma
2.2.2]{HH0}). Consider now the general case. The fact that the class
of local weak equivalences is closed under coproducts and filtered
colimits follows easily from the definition. Let $f:B\sr B'$ be a
morphism of bisimplicial presheaves such that the rows $f_i:B_i\sr
B_i'$ are local weak equivalences. Applying the wrapping functor (see
\cite[\S 4.2]{HH0}) to the simplicial objects in $\Delta^{op}PreShv(C)$
formed by the rows of $B$ and $B'$ we get a commutative square of the
form
$$
\begin{CD}
Wr_{col}B @>>> Wr_{col}B'\\
@VVV @VVV\\
B @>>> B'
\end{CD}
$$
Applying \cite[Proposition 4.2.8]{HH0} (see also \cite[Remark
4.3.9]{HH0}) we conclude that the vertical arrows are projective weak
equivalences and therefore define weak equivalences of the
corresponding diagonal objects. It remains to check that the diagonal
of the upper horizontal arrow is a local equivalence. This follows
easily from \cite[Lemma 4.2.3]{HH0},
\cite[Lemma 5.2.6 p. 126]{Hovey} and the fact (see \cite{J2}) that there
is a closed model structure on $\Delta^{op}PreShv(C)$ where local
equivalences are weak equivalences and all monomorphisms are
cofibrations.
\end{proof}
\begin{theorem}
\llabel{tmain} Let $P$ be a complete regular and bounded
cd-structure. Then the class of local weak equivalences in the
topology associated to $P$ coincides with the class
$$cl_{l}(W_{P})=cl_{\bdl}(W_{P}\cup W_{proj})$$ 
\end{theorem}
\begin{proof}
Lemma \ref{l7.31.3} together with Lemma \ref{es} implies that the
class of $t_P$-local weak equivalences contains
(\ref{eq7.31.2}). It remains to show that if $f:X\sr Y$ is a
$t_P$-local weak
equivalence then $f\in cl_{\bdl}(W_P\cup W_{proj})$. Let $N$ be
the set of morphisms $\emptyset'\sr U$ for all $U$ in $C$. Consider
the functor $Ex=Ex_{W_P,N}$ of \cite[Proposition 2.2.12]{HH0}. For any $f:X\sr Y$
we have a commutative square
$$
\begin{CD}
X @>f>> Y\\
@VVV @VVV\\
Ex(X) @>Ex(f)>> Ex(Y)
\end{CD}
$$
The vertical arrows are in $cl_{\bdl}(W_P)$ and thus they are local
weak equivalences. This implies that if $f$ is a local weak
equivalence then so is $Ex(f)$. On the other hand \cite[Proposition
2.2.12(2,3)]{HH0} together with Lemma \ref{l8.7.1} imply that the objects
$Ex(X)$ and $Ex(Y)$ are flasque. Thus by Lemma \ref{key} the morphism
$Ex(f)$ is a projective weak equivalence.
\end{proof}
Let $X\mapsto X_+$ be the functor from presheaves to pointed
presheaves left adjoint to the forgetful functor. For a class of
morphisms $E$ in $\Delta^{op}PreShv(C)$ denote by $E_+$ the class of
morphisms of the form $f_+:X_+\sr Y_+$ in the category of pointed
simplicial presheaves. The same reasoning implies the following
pointed version of Theorem \ref{tmain}.
\begin{theorem}
\llabel{tmainp}
Let $P$ be a complete regular and bounded cd-structure. Then the class
of pointed local weak equivalences in
the topology associated to $P$ coincides with the class 
$$cl_{l}(W_{P,+})=cl_{\bdl}(W_{P,+}\cup W_{proj})$$
\end{theorem}
Let $H(C_{t_P})$ be the homotopy category of simplicial (pre-)sheaves
on the site $(C, t_P)$ and $H_{\bullet}(C_{t_P})$ its pointed
analog. Combining Theorems \ref{tmain} and \ref{tmainp} with
\cite[Corollary 4.3.8]{HH0} we get the following description of these
categories.
\begin{cor}
\llabel{c7.31.4}
Under the assumptions of Theorem \ref{tmain} the functors
$$\Phi:\Delta^{op}C^{\coprod}\sr H(C_{t_P})$$
$$\Phi_{\bullet}:\Delta^{op}C^{\coprod}_+\sr H_{\bullet}(C_{t_P})$$
are localizations and one has
$$iso(\Phi)=cl_{\bdl}(W_P)$$
$$iso(\Phi_{\bullet})=cl_{\bdl}(W_{P_+})$$ 
\end{cor}

\subsection{The Brown-Gersten closed model structure}
Let $C$ be a category with a cd-structure $P$ and $Shv(C,t_P)$ the
category of sheaves of sets on $C$ with respect to the associated
topology. We say that a morphism in $\Delta^{op}Shv(C, t_P)$ is a
local (resp. projective) weak equivalence if it is a local
(resp. projective) weak equivalence as a morphism of simplicial
presheaves. If we define weak equivalences in $\Delta^{op}Shv(C,t_P)$
as the local weak equivalences with respect to $t_P$, cofibrations as
all monomorphisms and fibrations by the right lifting property we get
a closed model structure on $\Delta^{op}Shv(C,t_P)$ which we call the
Joyal-Jardine closed model structure (see \cite{Jo}, \cite{J1},
\cite{J2}). In this section we show that if $P$ is complete regular
and bounded then there is another closed model structure on
$\Delta^{op}Shv(C,t_P)$ with the same class of weak equivalences which
we call the Brown-Gersten closed model structure. When $C$ is the
category of open subsets of a topological space and $P$ is the
standard cd-structure we recover the closed model structure
constructed in \cite{KSB2}. The Brown-Gersten closed model structure
is often more convenient then the Joyal-Jardine one since it is
finitely generated. Our definitions and proofs follow closely the
ones given in \cite{KSB2}.

For a simplicial set we denote by the same letter the corresponding
constant simplicial sheaf on $C$. Recall that $\Lambda^{n,k}$ is the
simplicial subset of the boundary $\partial\Delta^{n}$ of the standard
simplex $\Delta^{n}$ which is the complement to the k-th face of
dimension $n-1$. Denote by $J_P$ the class of morphisms of the
following two types:
\begin{enumerate}
\item for any object $X$ of $C$ the morphisms $\Lambda^{n,k}\times
\rho(X)\sr \Delta^n\times \rho(X)$ 
\item for any distinguished square of the form (\ref{eq1}) the
morphisms 
$$\Delta^{n}\times\rho(A)\coprod_{\Lambda^{n,k}\times\rho(A)}\Lambda^{n,k}\times\rho(X)\sr
\Delta^{n}\times\rho(X)$$
\end{enumerate}
and by $I_P$ the class of morphisms of the following two types:
\begin{enumerate}
\item for any object $X$ of $C$ the morphisms $\partial\Delta^n\times
\rho(X)\sr \Delta^n\times \rho(X)$ 
\item for any distinguished square of the form (\ref{eq1}) the
morphisms 
$$\Delta^{n}\times\rho(A)
\coprod_{\partial\Delta^n\times\rho(A)}\partial\Delta^n\times\rho(X)\sr 
\Delta^{n}\times\rho(X)$$
\end{enumerate}
The following lemma is straightforward.
\begin{lemma}
\llabel{l8.2.1}
A morphism $f$ has the right lifting property with
respect to $J_P$ if and only if the maps $E(X)\sr B(X)$ are Kan
fibrations for all $X$ in $C$ and the maps 
$$E(X)\sr B(X)\times_{B(A)}E(A)$$ are Kan fibrations for all
distinguished squares $Q$ of the form (\ref{eq1}). 

A map $f$ has the
right lifting property with respect to $I_P$ if the maps $E(X)\sr
B(X)$ are trivial Kan fibrations for all $X$ in $C$ and the maps
$$E(X)\sr B(X)\times_{B(A)}E(A)$$ are Kan fibrations for all
distinguished squares $Q$ of the form (\ref{eq1}).
\end{lemma}
\begin{definition}
\llabel{fibrations} Let $P$ be a regular cd-structure. A morphism
$f:{E}\sr {B}$ in $\Delta^{op}Shv(C,t_P)$ is called a
Brown-Gersten fibration if it has the right lifting property with
respect to elements of $J_P$. 
\end{definition}
\begin{lemma}
\llabel{l7.31.1} Let $P$ be a regular cd-structure, then any
Brown-Gersten fibrant simplicial sheaf is flasque.
\end{lemma}
\begin{proof}
Let $F$ be a fibrant simplicial sheaf. Since $F$ is a sheaf we have
$F(0)=pt$ and therefore the first condition of Definition \ref{d4} is
satisfied. Lemma \ref{l8.2.1} implies that for any distinguished
square of the form (\ref{eq1}) the map of simplicial sets $F(X)\sr F(A)$
is a Kan fibration. On the other hand since $F$ is a sheaf and $P$ is
regular Proposition \ref{p8.2.2} shows that the square $F(Q)$ is a
pull-back square. This implies that it is a homotopy pull-back square.
\end{proof}
A morphism is called a Brown-Gersten cofibration of it has the left
lifting property with respect to all trivial Brown-Gersten fibrations
that is morphisms which are Brown-Gersten fibrations and (local) 
weak equivalences. The proof of the following lemma is parallel to the
proof of \cite[Lemma, p.274]{KSB2}.
\begin{lemma}
\llabel{trfib} Let $P$ be a complete, bounded and regular
cd-structure. Then a morphism $f:{E}\sr {B}$ is a trivial
Brown-Gersten fibration if and only if it has the right lifting
property with respect to elements of $I_P$.
\end{lemma}
\begin{proof}
Lemma \ref{l8.2.1} implies that any map which has the right lifting
property with respect to $I_P$ is a projective weak equivalence and
has the right lifting property with respect to $J_P$. This proofs the
``if'' part of the lemma.

\noindent
To prove the ``only if'' part we have to check that a morphism $f$
which has the right lifting property for elements of $J_P$ and which
is a Brown-Gersten fibrations has the right lifting property with
respect to $I_P$. Using again Lemma \ref{l8.2.1} one concludes that it
is sufficient to show that under our assumption $f$ is a projective
weak equivalence i.e. that for any $X$ in $C$ the map of simplicial
sets $E(X)\sr B(X)$ defined by $f$ is a weak equivalence. Since this
map is a fibration it is sufficient to check that its fiber over any
point $x_0\in B(X)$ is contractible. The point $x_0$ defines a map
$\rho(X)\sr B$ and the contractibility of the fiber is equivalent to
the condition that the pull-back $f_{x_0}:\rho(X)\times_B E\sr
\rho(X)$ of $f$ with respect to this map is a projective weak
equivalence. This map is a local weak equivalence in the
$t_P$-topology (if our site has enough points it is obvious; for the
general case see \cite[Th.1.12]{J2}) and has the right lifting
property with respect to $J_P$. Observe now that the morphism
$\rho(U)\sr pt$ has the right lifting property with respect to $J_P$
and therefore so does the morphism $\rho(X)\times_B E\sr pt$. We
conclude that $f_{x_0}$ is a local weak equivalence between two
fibrant objects and since, by Lemma \ref{l7.31.1}, Brown-Gersten
fibrant objects are flasque we conclude that $f_{x_0}$ is a projective
weak equivalence by Lemma \ref{key}.
\end{proof}
\begin{theorem}
\llabel{isclmod} Let $P$ be a complete, bounded and regular
cd-structure. Then the classes of local weak equivalences,
Brown-Gersten fibrations and Brown-Gersten cofibrations form a closed
model structure on $\Delta^{op}Shv(C,t_P)$.
\end{theorem}
\begin{proof}
Parallel to the proof of \cite[Th.2]{KSB2}.
\end{proof}
\begin{lemma}
\llabel{isamono} If $P$ is a regular cd-structure then any
Brown-Gersten cofibration is a monomorphism.
\end{lemma}
\begin{proof}
Let $f:{A}\sr {B}$ be a cofibration. Using the standard
technique we can find a decomposition $e\circ p$ of the morphism
${A}\sr pt$ such that $p$ is a trivial fibration and $e$ is
obtained from elements of $I_P$ by push-outs and infinite
compositions. Since $P$ is regular elements of $I_P$ are monomorphisms
and therefore $e$ is a monomorphism. The definition of cofibrations
implies that there is a morphism $g$ such that $g\circ f=e$ which
forces $f$ to be a monomorphism.
\end{proof}
\begin{proposition}
\llabel{iscel}
The Brown-Gersten closed model structure is finitely generated (in the
sense of \cite[Section 4]{HoveyStab}) and cellular (in the sense of
\cite{Hirs}) . 
\end{proposition}
\begin{proof}
Lemma \ref{compactness} implies easily that domains and codomains of
elements of $I_P$ and $J_P$ are finite. Together with our definition
of fibrations and Lemma \ref{trfib} it implies that the Brown-Gersten
closed model structure is finitely generated. A finitely generated
closed model structure is cellular if any cofibration is an effective
monomorphism. Therefore, the second statement follows from Lemma
\ref{isamono} and the fact that any monomorphism in the category of
sheaves is effective.
\end{proof}
\begin{proposition}
The Brown-Gersten closed model structure is both left and right
proper.
\end{proposition}
\begin{proof}
Consider a commutative square of the form
\begin{equation}
\llabel{asquare}
\begin{CD} 
{B} @>e'>> {Y}\\
@Vp'VV @VVpV\\
{A} @>e>> X
\end{CD}
\end{equation}
Assume first that it is a push-out square, $e'$ is a cofibration and
$p'$ is a local weak equivalence. By Lemma \ref{isamono} $e'$ is a
monomorphism and therefore $p$ is a local weak equivalence by
\cite[Prop. 1.4]{J3}.

Assume that (\ref{asquare}) is a pull-back square, $p$ is a fibration
and $e$ is a local weak equivalence. The definition of Brown-Gersten
fibrations implies in particular that for any object $U$ of $C$ the
morphism of simplicial sets ${Y}(U)\sr X(U)$ is a Kan
fibration. In particular it is a local fibration in the sense of
\cite{J2}. The same reasoning as in the proof of \cite[Prop. 1.4]{J3}
shows now that $p'$ is a local weak equivalence.
\end{proof}

\subsection{Cech morphisms}
\label{sec10}
Let $C$ be a category with fiber products. For any morphism $f:X\sr Y$
in $C$ denote by $\check{C}(f)$ the simplicial object with the terms 
$$\check{C}(f)_{n}=X^{n+1}_Y$$ and faces and degeneracy morphisms
given by partial projections and diagonals respectively. The
projections $X^{i}_Y\sr Y$ define a morphism $\check{C}(f)\sr Y$ which
we denote by $\eta(f)$. For a class of morphisms $A$ in $C$ we denote
by $\eta(A)$ the class of morphisms of the form $\eta(f)$ for $f\in
A$. 

\noindent
For a cd-structure $P$ on a category $C$ denote by $Cov_f(t_P)$ the
class of morphisms in $C^{{\coprod}_{<\infty}}$ of the form $f=\coprod
f_i$ for all finite $t_P$-coverings $\{f_i:X_i\sr X\}$ in $C$. Note
that since the empty family is a covering of the initial object the
morphism $\emptyset'\sr \emptyset$ belongs to $Cov_f(t_P)$. If $C$ has
fiber products then so does $C^{\coprod_{<\infty}}$ and therefore we
may consider the class of morphisms $\eta(Cof_f(t_P))$. The goal of
this section is to show that, for a good enough cd-structure $P$ on a
category $C$ with fiber products, the class
$cl_{\Delta,\coprod_{<\infty}}(W_P)$ coincides with the class
$cl_{\Delta,\coprod_{<\infty}}(\eta(Cov_t(t_P)))$. 
\begin{lemma}
\llabel{section}
Let $f:X\sr Y$ be a morphism in $C$ which has a
section. Then $\eta(f)$ is a homotopy equivalence.
\end{lemma}
\begin{proof}
If $f$ has a section
$Y\sr X$ it defined a map $Y\sr \check{C}(f)$ in the obvious way. The
composition $Y\sr \check{C}(f)\sr Y$ is the identity. On the other
hand, for any simplicial object $K$ over $Y$ one has
$$Hom_{Y}(K,\check{C}(f))=Hom_{Y}(K_0,X)$$
where $K_0$ is the object of 0-simplexes of $K$. In particular any two
morphisms with values in $\check{C}(f)$ over $Y$ are homotopic. This
implies that the composition $\check{C}(f)\sr Y\sr \check{C}(f)$ is
homotopic to identity. 
\end{proof}
\begin{lemma}
\llabel{b130}
Let $C$ be a category with fiber products and $Q$ a fiber square of
the form (\ref{eq1}) such that both morphisms $A\sr X$ and $Y\sr X$
are monomorphisms. Denote by $\pi_Q$ the morphism $Y\coprod A\sr
X$. Then one has
$$p_Q\in cl_{\Delta,\coprod_{<\infty}}(\eta(\pi_Q))$$
\end{lemma}
\begin{proof}
Consider the fiber square
\begin{equation}
\llabel{b13sq}
\begin{CD}
K_Q\times_X \check{C}(\pi_Q)@>>>\check{C}(\pi_Q)\\
@VVV @VVV\\
K_Q @>p_Q>> X
\end{CD}
\end{equation}
Let us show that the upper horizontal and the left vertical arrows belong to
$cl_{\Delta, \coprod_{<\infty}}(\emptyset)$. Thinking of the fiber
product as of the diagonal of the corresponding bisimplicial object we
see that the left vertical arrow is contained in the $\Delta$-closure
of the family of morphisms of the form
$$(Y\coprod A)^n\times\check{C}(\pi_Q)\sr (Y\coprod A)^n$$
where the products are taken over $X$.
Since 
$$(Y\coprod A)^n_X=\coprod_{i+j=n} Y^i\times_X A^j$$
it is contained in the $(\Delta,\coprod_{<\infty})$-closure of
morphisms of the form 
$$Y^i\times_X A^j\times_X\check{C}(\pi_Q)\sr Y^i\times_X A^j$$
which are isomorphic to morphisms of the form $\eta(\pi_Q\times_X
Id_{Y^i\times_X A^j})$. Since $i+j=n>0$ the morphism $\pi_Q\times_X
Id_{Y^i\times_X A^j}$ has a section. Lemma \ref{section} implies now
that $\eta(\pi_Q\times_X
Id_{Y^i\times_X A^j})$ are homotopy equivalences and therefore the
left vertical arrow is in $cl_{\Delta,
\coprod_{<\infty}}(\emptyset)$. 

\noindent
For the upper horizontal arrow the same argument shows that it belongs
to the $(\Delta,\coprod_{<\infty})$-closure of morphisms of the form
$p_{Q'}$ where $Q'=Q\times_{X}Y^i\times_X A^j$ and $i+j>0$. Since both
$A\sr X$ and $Y\sr X$ are monomorphisms each of these squares has the
property that its vertical or horizontal sides are isomorphism. For
such squares the morphisms $p_Q$ are homotopy equivalences by
\cite[Lemma 2.1.5]{HH0}. 
\end{proof}
\begin{lemma}
\llabel{b131}
Let $C$ be a category with fiber products and $P$ a regular
cd-structure on $C$. Then for any distinguished square $Q$ one has
$$p_Q\in cl_{\Delta,\coprod_{<\infty}}(\eta(Cov(t_P))).$$
\end{lemma}
\begin{proof}
Let $\pi_Q:Y\coprod A\sr X$ be the element of $Cov(t_P)$ defined by
$Q$. Consider again a square of the form (\ref{b13sq}). The same
argument as in the proof of Lemma \ref{b130} shows that the left
vertical arrow is in $cl_{\Delta, \coprod_{<\infty}}(\emptyset)$ and
the upper horizontal arrow belongs to the
$(\Delta,\coprod_{<\infty})$-closure of morphisms of the form $p_{Q'}$
where $Q'=Q\times_{X}Y^i\times_X A^j$ and $i+j>0$. For since $A\sr X$
is a monomorphism by the definition of a regular cd-structure for
$j>0$ the horizontal arrows of $Q'$ are isomorphisms and the
corresponding morphism is a homotopy equivalence by \cite[Lemma
2.1.5]{HH0}. Consider the squares $Q_i=Q\times_{X}Y^i$ for
$i>0$. Since $A\sr X$ is a monomorphism we have $B\times_X Y=B\times_X
B$ and therefore the the square $Q_1$ is isomorphic to the lower
square of the following diagram
\begin{equation}
\begin{CD}
B@>>>Y\\
@VVV @VVV\\
B\times_X B@>>>Y\times_X Y\\
@VVV @VVV\\
B@>>>Y
\end{CD}
\end{equation}
where the upper vertical arrows are the diagonals. Denote the upper
square of this diagram by $Q'$. By \cite[Lemma 2.1.10]{HH0} we
conclude that it is sufficient to show that one has
\begin{equation}
\llabel{b13in}
p_{Q'\times_X Y^i}\in
cl_{\Delta,\coprod_{<\infty}}(\eta(Cov(t_P)))
\end{equation}
for $i\ge 0$. Both morphisms in $Q'$ are monomorphisms and $B\times_X
B\coprod Y\sr Y\times_X Y$ is in $Cov_{t_P}$ by definition of a
regular cd-structure. Therefore the same holds for $Q'\times_X Y^i$
and (\ref{b13in}) follows from Lemma \ref{b130}.
\end{proof}
\begin{proposition}
\llabel{cmain}
Let $A$ be a class of morphisms in $C$ which is closed
under pull-backs. Then one has:
\begin{enumerate}
\item if $f,g$ is a pair of morphisms such that $gf$ is defined and 
$\eta(gf)\in cl_{\Delta}(\eta(A))$ then $\eta(g)\in cl_{\Delta}(\eta(A))$ 
\item if $f,g$ is a pair of morphisms such that $gf$ is defined and $\eta(f),
\eta(g)\in 
cl_{\Delta}(\eta(A))$ then $\eta(gf)\in cl_{\Delta}(\eta(A))$.
\end{enumerate}
\end{proposition}
\begin{proof}
The first statement follows from Lemma \ref{trrr}. The second one from
Lemma \ref{comp}.
\end{proof}
\begin{lemma}
\llabel{gen0}
Let $p:A\sr X$, $q:B\sr X$ be two morphisms. Then one has
\begin{equation}
\label{incl}
\eta(p)\in cl_{\Delta}(\{\eta({p}\times_{X} Id_{{B}^m}),
\eta({q}\times_{X} Id_{A^{n}})\}_{m>0,n\ge 0})
\end{equation}
\end{lemma}
\begin{proof}
Consider the bisimplicial object $\check{C}(p,q)$ build on $p$ and
$q$. We have a commutative diagram
\begin{equation}
\llabel{cpq}
\begin{CD}
\Delta\check{C}(p,q)@>>>\check{C}(q)\\
@VVV @VVV\\
\check{C}(p)@>>>X
\end{CD}
\end{equation}
The upper horizontal arrow belongs to
$cl_{\Delta}(\{\eta({p}\times_{X} Id_{B^m})\}_{m>0})$ and the
left vertical one to 
$cl_{\Delta}(\{\eta({q}\times_{X} Id_{A^n})\}_{n>0})$. This
implies (\ref{incl}) by the ``2 out of 3'' property of the
$\Delta$-closed classes.
\end{proof}
\begin{lemma}
\llabel{trrr}
Let $X\stackrel{f}{\sr}Y\stackrel{g}{\sr}Z$ is a composable pair of
morphisms. Then one has
\begin{equation}
\label{incl2}
\eta(g)\in
cl_{\Delta}(\{\eta(gf\times_{Z}Id_{Y^n})\}_{n\ge 0})
\end{equation}
\end{lemma}
\begin{proof}
Applying Lemma \ref{gen0} to the pair $p=g, q=gf$ we get
$$\eta(g)\in cl_{\Delta}(\{\eta({g}\times_{Z} Id_{{X}^m}),
\eta({(gf)}\times_{Z} Id_{Y^{n}})\}_{m>0,n\ge 0})$$
Since morphisms ${g}\times_{Z} Id_{{X}^m}$ for $m>0$ have sections
Lemma \ref{section} implies (\ref{incl2}).
\end{proof}
\begin{lemma}
\llabel{comp}
Let $X\stackrel{f}{\sr}Y\stackrel{g}{\sr}Z$ is a composable pair of
morphisms. Then one has
\begin{equation}
\label{incl3}
\eta(gf)\in cl_{\Delta}(\{\eta(g\times
Id_{X^n}),\eta(f\times Id_{X^n}\times Id_{Y^l})\}_{n,l\ge
0})
\end{equation}
where all the products are taken over $Z$.
\end{lemma}
\begin{proof}
All the products in the proof are over $Z$ unless the opposite is
specified. Let us show first that for $l\ge 0$ one has
\begin{equation}
\label{incl4}
\eta((gf)\times
Id_{Y^{l+1}})\in cl_{\Delta}(\eta(f\times
Id_{X^n}\times Id_{Y^l})\}_{n,l\ge 0})
\end{equation}
For that we apply Lemma \ref{gen0} to $p=(gf)\times Id_{Y^{l+1}}$
and $q=(Y^l\times X\stackrel{Id\times f}{\sr} Y^{l}\times Y)$. We have
\begin{equation}
\llabel{morph1}
((gf)\times Id_{Y^{l+1}})\times_{Y^{l+1}}
Id_{{(Y^{l}\times X)}^m}=(gf)\times Id_{{(Y^{l}\times X)}^m}
\end{equation}
and
$${q}\times_{Y^{l+1}} Id_{({X\times
Y^{l+1})}^{n}}=q\times Id_{X^n}=f\times Id_{X^n}\times Id_{Y^l}$$
and therefore 
$$\eta((gf)\times Id_{Y^{l+1}})\in
cl_{\Delta}(\{\eta((gf)\times Id_{{(Y^{l}\times X)}^m}),
\eta(f\times Id_{X^n}\times Id_{Y^l})\}_{m>0,n\ge 0})$$
For $m>0$ the morphisms (\ref{morph1}) have sections and therefore Lemma
\ref{section} implies (\ref{incl4}). Apply now Lemma \ref{gen0} to
morphisms $p=gf$ and $q=g$. We get
$$\eta(gf)\in cl_{\Delta}(\{\eta({gf}\times Id_{{Y}^m}),
\eta({g}\times Id_{X^{n}})\}_{m>0,n\ge 0})$$
which together with (\ref{incl4}) implies (\ref{incl3}).
\end{proof}
\begin{proposition}
\llabel{complcase} Let $P$ be a cd-structure on a category $C$ with
fiber products such that for any distinguished square $Q$ and a
morphism $X'\sr X$ the square $Q\times_X X'$ is distinguished. Then
for any covering $\{f_i:X_i\sr X\}$ in the associated topology one has
$$\eta(Cov(t_P))\subset cl_{\Delta,\coprod_{<\infty}}(W_P)$$
\end{proposition}
\begin{proof}
Our assumption on $P$ implies in particular that it is complete. The
class $Cov(t_P)$ is then the smallest class which contains the
morphism $\emptyset'\sr \emptyset$, morphisms of the form
$\pi_Q:Y\coprod A\sr X$ for distinguished squares $Q$ of the form
(\ref{eq1}) and is closed under finite coproducts and under operations
described in parts (1) and (2) of Proposition \ref{cmain}. This
implies that it is sufficient to prove that for any $Q$ in $P$ one has
$$\eta(\pi_Q)\in cl_{\Delta,\coprod_{<\infty}}(\{p_{Q'}\}_{Q'\in
P})$$
This follows from the diagram (\ref{b13sq}) in the same way as in the
proof of Lemma \ref{b130}. 
\end{proof}
Combining Propositions \ref{b131} and \ref{complcase} we get the
following result.
\begin{proposition}
\llabel{m94} Let $C$ be a category with fiber products and $P$ be a
regular cd-structure such that for any distinguished distinguished
square $Q$ and a morphism $X'\sr X$ the square $Q\times_X X'$ is
distinguished. Then one has
$$cl_{\Delta,\coprod_{<\infty}}(\eta(Cov_f(t_P)))=
cl_{\Delta,\coprod_{<\infty}}(W_P)$$   
\end{proposition}
Let $Cov(t_P)$ be the class of morphisms in $C^{\coprod}$ of the form
$\coprod f_i$ for all (as opposed to just finite) coverings of the
form $\{f_i:U_i\sr U\}$ in $C$.
\begin{cor}
\llabel{c8.3.11} Under the assumptions of the
proposition one has
$$cl_{\bdl}(W_P)=cl_{\bdl}(\eta(Cov(t_P)))$$
\end{cor}
\begin{proof}
We have only to show that for any covering $\{f_i\}$ the morphism
$\eta(f)$ is in $cl_{\bdl}(W_P)$. This follows from the proposition
and Lemma \ref{trrr} since any covering has a finite refinement.
\end{proof}

\subsection{Reformulations for categories with finite coproducts}
Let $C$ be a category with an initial object $0$ and finite
coproducts. Recall (see \cite{HH0}) that we denote by $R(C)$ the
category of radditive functors on $C$ and by $[C]$ the full
subcategory in $R(C)$ which consists of coproducts of representable
functors. Let $P$ be a cd-structure on $C$. In this section we show
how, under some conditions on $C$ and $P$, the results of the previous
sections can be reformulated in terms of $R(C)$ instead of $PreShv(C)$
and $[C]$ instead of $C^{\coprod}$.
\begin{definition}
\llabel{chunky}
A small category $C$ is called chunky if it has an initial object $0$
and finite coproducts and the following conditions hold:
\begin{enumerate}
\item any morphism with values in $0$ is an isomorphism
\item any square of the form
\begin{equation}
\llabel{eq8.3.0}
\begin{CD}
0 @>>> X\\
@VVV @VVV\\
Y @>>> X\coprod Y
\end{CD}
\end{equation}
is pull-back
\item for any objects $X$, $Y$, $Z$ and a morphism $Z\sr X\coprod Y$
the fiber products $Z\times_{X\coprod Y} X$ and $Z\times_{X\coprod Y}
Y$ exist and the morphism
$$(Z\times_{X\coprod Y} X)\coprod (Z\times_{X\coprod Y}
Y)\sr Z$$
is an isomorphism.
\end{enumerate}
\end{definition}
For a category $C$ with finite coproducts and an initial object denote
by $P_{add}$ the cd-structure which consists of squares of the form
(\ref{eq8.3.0}). 
\begin{lemma}
\llabel{l8.3.1}
Let $C$ be a chunky category. Then $P_{add}$ is complete, regular and
bounded.
\end{lemma}
\begin{proof}
The fact that $P_{add}$ is complete follows from the first and the
third conditions of Definition \ref{chunky} and Lemma
\ref{completepb}. To check that it is regular it is clearly sufficient
to verify that for any $X$ and $Y$ the morphism $X\sr X\coprod Y$ is a
monomorphism. This follows easily from our conditions.
To show that $P_{ad}$ is bounded consider the density structure $D$
such that $D_0(U)$ consists of all morphisms with values in $U$ and
$D_n(U)$ for $n\ge 1$ consists of all isomorphisms with values in
$U$. One verifies easily that any square of the form (\ref{eq8.3.0})
is reducing for this density structure.
\end{proof}
\begin{lemma}
\llabel{l8.3.2}
Let $C$ be a chunky category. Then a presheaf on $C$ is a radditive
functor if and only if it is a sheaf in the topology $t_{add}$
associated to $P_{add}$.
\end{lemma}
\begin{proof}
Follows from Lemma \ref{l8.3.1} and Corollary \ref{regcompl}.
\end{proof}
\begin{lemma}
\llabel{l8.3.3.1}
Let $C$ be a chunky category, $F$ a presheaf on $C$ and $r(F)$ the
associated radditive functor. Then the morphism $F\sr r(F)$ belongs to
$cl_{\bdl}(W_{P_{add}}\cup W_{proj})$.
\end{lemma}
\begin{proof}
Lemma \ref{l8.3.2} implies that any morphism of the form $F\sr r(F)$
is a local equivalence in $t_{add}$. Our result follows now from
Theorem \ref{tmain} and Lemma \ref{l8.3.1}.
\end{proof}
For an object $U$ in $C$ denote by $D/U$ the category of the following
form
\begin{enumerate}
\item an object of $D/U$ is a sequence of morphisms $(u_1:U_1\sr
U,\dots, u_n:U_n\sr U)$ in $C$ such that $U_n\ne 0$ and $\coprod u_n$
is an isomorphism
\item a morphism from $(u_i:U_i\sr
U)_{\{i\in I\}}$ to $(v_j:V_j\sr
U)_{\{j\in J\}}$ is a surjection $\alpha:I\sr J$ together with a
collection of morphisms $U_i\sr V_{\alpha(i)}$ over $U$ such that
for every $j\in J$ the morphism 
$$\coprod_{i\in \alpha^{-1}(j)} U_i\sr V_j$$
is an isomorphism.
\end{enumerate} 
\begin{lemma}
\llabel{l8.3.3} Let $C$ be a chunky category and $U\ne 0$ an object in
$C$. Then the category $D/U$ is filtered and for any presheaf $F$ on
$C$ one has
\begin{equation}
\llabel{eq8.3.4}
r(F)(U)=colim_{(U_1,\dots, U_n)\in D/U} F(U_1)\times\dots\times F(U_n)
\end{equation}
\end{lemma}
\begin{proof}
The fact that coprojections in $C$ are monomorphisms, squares of the
form (\ref{eq8.3.0}) are pull back and $Hom(U,0)=\emptyset$ for $U\ne
0$ imply that the category $D/U$ is a partially ordered set
(i.e. there is at most one morphism between two objects). The third
condition of Definition \ref{chunky} implies that for any two objects
in $D/U$ there exists a third one which maps to both. Together, these
two facts imply that $D/U$ is filtered. 

Using again the third condition of Definition \ref{chunky} one shows
that the right hand side of (\ref{eq8.3.4}) defines a functor $F'$ on
$C$. One verifies easily that it is radditive and that for any
radditive $G$ the obvious map
$$Hom(F',G)\sr Hom(F,G)$$
is a bijection.
\end{proof}
\begin{lemma}
\llabel{l8.3.5}
Let $C$ be a chunky category. The the functor $F\mapsto r(F)$ takes
projective weak equivalences of simplicial presheaves to projective
weak equivalences.
\end{lemma}
\begin{proof}
It follows from Lemma \ref{l8.3.3} because the class of weak
equivalences of simplicial sets is closed under finite products and
filtered colimits.
\end{proof}
For a cd-structure $P$ on a category $C$ denote by $H(C,P)$ the
localization of $\Delta^{op}PreShv(C)$ with respect to
$cl_{\bdl}(W_P\cup W_{proj})$. 
\begin{proposition}
\llabel{p8.3.6}
Let $C$ be a chunky category and $P$ a cd-structure on $C$ which
contains $P_{add}$.  Then the functor 
$$\Delta^{op}PreShv(C)\sr H(C,P)$$
factors through a functor
\begin{equation}
\llabel{eq8.3.7}
\Phi:\Delta^{op}R(C)\sr H(C,P)
\end{equation}
which is a localization and $iso(\Phi)=cl_{\bdl}(r(W_P)\cup W_{proj})$.
\end{proposition}
\begin{proof}
The functor $r:\Delta^{op}PreShv(C)\sr \Delta^{op}R(C)$ is a left
adjoint to a full embedding and therefore it is a localization with
respect to morphisms of the form $F\sr r(F)$. These morphisms map to
isomorphisms in $H(C,P)$ by Lemma \ref{l8.3.3.1} and our assumption
that $P$ contains $P_{add}$. This implies that $\Phi$ is defined and
that it is a localization. The inclusion 
$$cl_{\bdl}(r(W_P)\cup W_{proj})\subset iso(\Phi)$$
is obvious. To prove the opposite inclusion consider an element
$f:X\sr Y$ of $iso(\Phi)$. As a morphism in $\Delta^{op}PreShv(C)$ it
belongs to  $cl_{\bdl}(W_P\cup W_{proj})$ by \cite{HH0}. Since $r$
commutes with colimits we have
$$r(cl_{\bdl}(W_{P}\cup W_{proj}))\subset cl_{\bdl}(r(W_P)\cup W_{proj})$$
and since $r(f)=f$ this implies the required result.
\end{proof}
\begin{lemma}
\llabel{oldchunky}
Let $C$ be a chunky category. Then any simplicial radditive functor on
$C$ is vertically clant. In particular, any chunky category is grainy.
\end{lemma}
%
%
\begin{lemma}
\llabel{l8.3.8}
Let $C$ and $P$ be as in Proposition \ref{p8.3.6}. Then the functor
$$\Phi:\Delta^{op}[C]\sr H(C,P)$$
is a localization and $iso(\Phi)=cl_{\bdl}(W_P)$.
\end{lemma}
\begin{proof}
This functor is a composition of the form
$$\Delta^{op}[C]\sr \Delta^{op}R(C)\sr H(C)\sr H(C,P)$$
The composition of the second and the third arrow is a localization by
Proposition \ref{p8.3.6}, the second arrow is a localization with
respect to $W_{proj}$ by definition. Therefore, the third arrow is a
localization. On the other hand, the composition of the first and the
second arrows is a localization by \cite[Proposition 3.4.9]{HH0} and
therefore $\Phi$ is a localization. Let $f$ be a morphism in
$\Delta^{op}[C]$ which maps to an isomorphism in $H(C,P)$. By
Proposition \ref{p8.3.6} its image in $\Delta^{op}R(C)$ is in
$cl_{\bdl}(W_P\cup W_{proj})$. Therefore, by Lemma \ref{oldchunky} and
\cite[Proposition 4.1.15]{HH0} we have $f\in cl_{\bdl}(W_P)$.
\end{proof}
\begin{lemma}
\llabel{l8.3.9} 
Let $C$ and $P$ be as in Proposition
\ref{p8.3.6}. Assume in addition that $C$ has fiber products, $P$ is
regular and that the pull-backs of distinguished squares are
distinguished. Then one has
\begin{equation}\llabel{eq8.3.10}
cl_{\Delta,\coprod_{<\infty}}(\eta(Cov_{f}(t_P)))=
cl_{\Delta,\coprod_{<\infty}}(W_P) 
\end{equation}
where the corresponding classes are considered in $\Delta^{op}C$.
\end{lemma}
\begin{proof}
By Proposition \ref{m94} an equality of the form (\ref{eq8.3.10})
holds in $\Delta^{op}C^{\coprod_{<\infty}}$. The functor
$C^{\coprod_{<\infty}}\sr C$ commutes with finite coproducts and
therefore takes elements of the class $W_P$ in
$\Delta^{op}C^{\coprod_{<\infty}}$ to elements of the analogous class
in $\Delta^{op}C$.  The third condition of Definition \ref{chunky}
implies that this functor also commutes with
fiber products and therefore takes elements of
$\eta(Cov_{f}(t_P))$ in $\Delta^{op}C^{\coprod_{<\infty}}$ to
elements of $\eta(Cov_{f}(t_P))$ in $\Delta^{op}C$ which implies the
statement of the lemma.
\end{proof}

\bibliography{alggeom}
\bibliographystyle{plain}
\end{document}